%% file: KCL_review_article.tex
\documentclass[12pt] {amsart}

\date{}
\usepackage{graphicx}
\input psfigfx
\usepackage{url}
\numberwithin{equation}{section}

\usepackage{xcolor}
\usepackage[mathscr]{euscript}

\newcommand{\red}[1]{\textcolor{red}{#1}}
\newcommand{\blue}[1]{\textcolor{blue}{#1}}
\newcommand{\green}[1]{\textcolor{green}{#1}}

\renewcommand{\baselinestretch}{1.2}

\newtheorem{theo}{Theorem}[section]

\newtheorem{defn}[theo]{Definition}

\font\bbb=msbm10 scaled\magstep1

\def\bold#1{\mbox{\boldmath $#1$}}

\newcommand{\RR}{\mbox{\bbb R}}
\newcommand{\uu}[1]{\bold{#1}}

\begin{document}
\newcommand{\bea}{\begin{eqnarray}}
\newcommand{\eea}{\end{eqnarray}}

\title[Conserved Variable and Flux in KCL]{System of Kinematical Conservation Laws (KCL) \\ A Review Article}

\author[Arun]{K. R. Arun}
\address{School of Mathematics, Indian Institute of Science Education and Research Thiruvananthapuram, Thiruvananthapuram 695551, India}
\email{arun@iisertvm.ac.in}

\author[Prasad]{Phoolan Prasad}
\address{Retired, Department of Mathematics, Indian Institute of Science, Bangalore 560012, India}
\email{phoolan.prasad@gmail.com}
\urladdr{//math.iisc.ernet.in/~prasad/}

\date{\today}

\keywords{System of Kinematical Conservation Laws (KCL), Ray Theory,  Conservation Laws, Propagation of Curves and Surfaces}

\begin{abstract}
In a wide range of physical phenomena, we find propagating surfaces $\Omega_t$ which need mathematical treatment.  In this article, we review the theory of the system of kinematical conservation laws (KCL), which govern the evolution of these surfaces and have been developed by the second author and his collaborators. KCL are the most general equations in conservation form, governing the evolution of $\Omega_t$ with physically realistic singularities. A very special type of singularity is a kink, which is a point on $\Omega_t$ when $\Omega_t$ is a curve in $\RR^2$ and is a curve on $\Omega_t$ when $\Omega_t$ is a surface in $\RR^3$. Across a kink the normal $\uu{n}$ to $\Omega_t$ and the normal velocity $m$ on $\Omega_t$ are  discontinuous. The main aim of this article is to identify density of the conserved variable and the flux for the KCL which we did not do earlier. The presentation of this article is like that in a popular article, which which aims at non-experts in the field. 

\end{abstract}

\maketitle

\section{Intoduction}

Consider a moving curve in two space dimensions (which we denote briefly by $2$-D), i.e.\ $(x,y)$-plane or a moving surface in three space dimensions (denote it by $3$-D), i.e.\ $(x,y,z)$-space. We denote this curve or surface\footnote{The results can be extended easily to a surface in space of arbitrary dimensions \cite{pp-16b, pp-book-18}} at time $t$ by $\Omega_t$. The movement of $\Omega_t$ is determined by movement of its points according to some law depending on the medium in which $\Omega_t$ evolves. The path of a point on $\Omega_t$ is called a ray. Once the so called ray velocity, $\uu \chi = (\chi_1, \chi_2)$ say, is known, (see the equation (\ref{5.2}) below) we can easily write an eikonal equation (see(\ref{5.5})) of the surface, which is a Hamilton-Jacobi PDE.
There is now a standard theory to study viscosity solution of this first order PDE, where the folded trailing part of $\Omega_t$ as seen in Figure 1, if formed by ray theory, disappears. However, when the rays have built in them effect of genuine nonlinearity of the equations governing a medium in which $\Omega_t$ propagates, a new type of singularity, called kink appears on $\Omega_t$ \cite{pp-book} (as seen in Figure 4 where two kinks appear and are shown by dots). Both the ray theory and the level set theory to solve the eikonal equation are inadequate to study formation and evolution of a kink.  

\begin{figure}
\centering
\includegraphics[height=0.450\textheight]{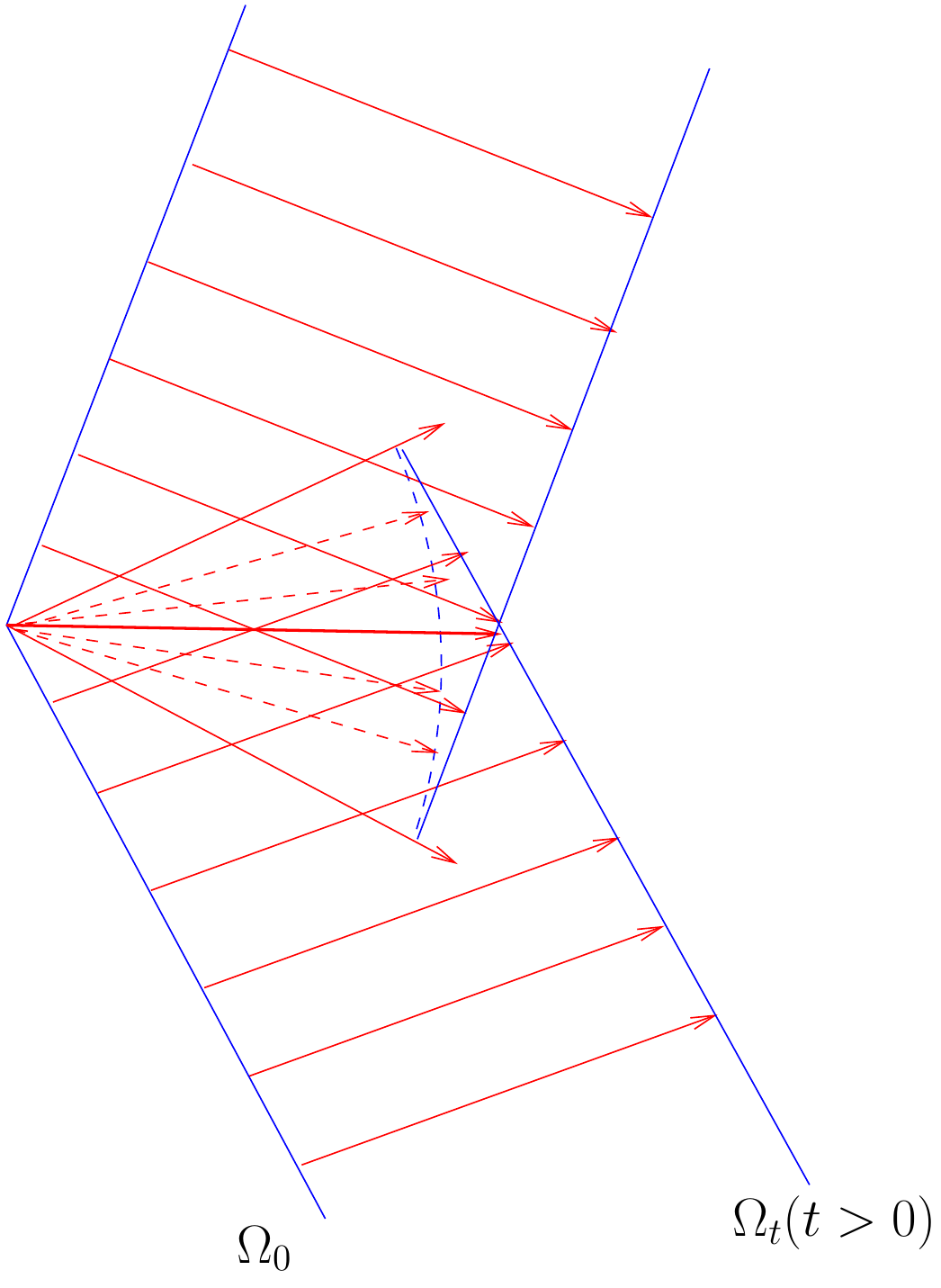}
\caption{Caustic in a Linear Theory, when the ray velocity is constant and normal to $\Omega_t$. The initial wavefront is formed by two intersecting straight lines forming a wedge concave to the direction of propagation. \textcolor{blue}{.....} part of the curve produced by the corner O  and determined by Huygen's method. This is an arc of a circle with center at O. \textcolor{blue}{-----} part of the curve produced by the smooth part of the initial position of the curve by ray theory.}    
\end{figure}
  
Distinguishing feature of a kink is in appearance of a discontinuity in the normal direction $\uu{n}$ and the velocity $m$ of $\Omega_t$.  Examples of $\Omega_t$ with discontinuities in $\uu{n}$ and $m$ across a $2$ dimensional surface on $\Omega_t$ are plenty. They were observed in experimental results \cite{stu-kul-76} and it became our aim to capture the experimental results mathematically. When discontinuities in $\uu n$ and $m$ appear on $\Omega_t$, the governing PDEs of $\Omega_t$, i.e.\, the eikonal equation and ray theory breakdown, then we need to go to the more basic formulation for the evolution of a curve or surface, which is KCL. \textbf{By KCL we mean ``a system of kinematical conservation laws"}. Conservation laws are formulation in terms of integrals, which remain valid even if singularities (like discontinuities in functions) on $\Omega_t$ appear. $2$-D KCL for evolution of a curve was formulated (rather first discovered by chance and then formulated) in 1992 by K. W. Morton, Phoolan Prasad and Renuka Ravindran and has been extensively used by Prasad and his collaborators to find out many new properties of weakly nonlinear wave fronts and shock fronts. The results have been published in many papers in important journals like Journal of Fluid Mechanics, IMA Journal of Applied Mathematics, Proceedings of Indian Academy of Sciences : Mathematical sciences, and many international conferences. $3$-KCL was first formulated in 1994 by Mike Giles, Phoolan Prasad and Renuka Ravindran but it was completed much later by K. R. Arun and Prasad. The papers and applications of $3$-D KCl were were published since 2009 in journals like Wave Motion, SIAM J. Sci. Comput. and international conferences. We do not wish to  give all these references but most of them are available in \cite{pp-book-18}. But we wish to emphasize that all these interesting and physically realistic results were obtained without identifying density of the conserved variable and flux for the KCL. We wish to identify these and complete the KCL theory in this article.

For readers not familiar with the theory of conservation laws, we start with a simple example, namely conservation of mass in fluid dynamics in the next section.

\section{Conservation of mass in fluid mechanics:} Let us consider fluid flow in an one-dimensional pipe along $x$-axis and consider two sections at $x_\ell$ and $x_r$. Let $\rho(x,t)$ represent the mass density and $u(x,t)$ the fluid velocity, which can be positive or negative. The total fluid contained between the two sections at time $t$ is $\int_{x_\ell}^{x_r}\rho(\xi,t)d\xi$.  The flux (per unit time) of mass at $x_\ell$ is $\rho(x_\ell)u(x_\ell)$ and similarly at $x_r$ is $\rho(x_r)u(x_r)$. When no fluid is created (for example by chemical reaction) balance between the rate of change of mass between these two sections and the fluxes at the two ends is expressed by the conservation law.
\begin{equation}
\label{2.1}
\frac{d}{dt}\int_{x_\ell}^{x_r}\rho(x,t)dx = \rho(x_\ell,t)u(x_\ell,t)-\rho(x_r,t)u(x_r,t).
\end{equation}
When the mass density and fluid velocity are differentiable, this equation immediately leads to (see the next section)
\begin{equation}
\label{2.2}
\rho_t + (\rho u)_x =0
\end{equation}
which gives the Euler's equation of continuity
\begin{equation}
\label{2.3}
\rho_t +  u \rho_x + \rho u_x = 0.
\end{equation}
As explained in the next section, we call (\ref{2.2}) the conservation law of mass but by this we mean the (\ref{2.1}). (\ref{2.3}) is the differential form of the conservation law of mass. Note that (\ref{2.1}) is incomplete since there are two state variables $\rho(x,t)$ and $u(x,t)$ and there is only one conservation law. 

\section{Scalar conservation law and curves of discontinuity in one-space dimension}
Let $H(u(x,t))$ the conserved variable associated with a state variable\footnote{By state variables we mean the variables which are required to identify the state of the system governed by the the equation or system of equations} $u(x,t)$  and $F(u(x,t))$ be the flux. We assume $H(u)$ and $F(u)$ to be smooth functions\footnote{We shall use the word `smooth' quite often in this article as we have done in the caption of the Figure 1 for a curve. By a smooth function  we mean a function that is as many times differentiable as required. By smooth a curve or smooth evolution of a curve we mean that the functions appearing in their description are smooth.} of $u$ . A scalar or single conservation law is 
\begin{equation}
{\frac{d}{dt}} \int\limits_{x_\ell}^{x_r} H(u(x, t))
dx = F(u(x_\ell, t)) - F(u(x_r, t)), ~~~\ x_\ell, x_r \ \hbox{fixed}.
 \label{3.1}
\end{equation} 
When $u$ is differentiable, we can show that\footnote{The equation (\ref{3.1}) can be written as
\begin{equation*}
\int\limits_{x_\ell}^{x_r} H_t(u(x, t))
dx = -\int_{x_\ell}^{x_r}F_x(u(x,t))dx.
\end{equation*} 
Since the interval $(x_\ell,x_r)$ is arbitrary (\ref{3.2}) follows. Note that $F(u(x_\ell, t)) - F(u(x_r, t)) = \int_{x_\ell}^{x_r}F_x(u(x,t))dx$ is the Gauss divergence theorem in one dimension. } (\ref{3.1}) reduces to 
\begin{equation}
\frac{\partial H(u)}{\partial t} + \frac{\partial F(u)}{\partial x} = 0.
\label{3.2}
\end{equation}
 \begin{defn}
We call (\ref{3.2}) a \textbf{conservation law}\footnote{Whether $u$ is differentiable or not.}, by which we mean the integral balance (\ref{3.1}). 
\end{defn}

We can show that {\it for a $\mathcal{C}^1$ function $u$, the above conservation
law is equivalent to the PDE}
\begin{equation}  
\label{3.3}
\dot H (u)u_t+\dot F(u)u_x=0  ~~\textit{or}~~ u_t + \frac{\dot F(u)} {\dot H(u)}u_x=0 ~~\textit{assuming}~~ \dot H(u)\neq 0,
\end{equation}

\noindent where dot above a quantity represents its derivative.

It is well known\footnote{See any standard book on PDE for example \cite{pp-rr-book} which is freely available on a cite} that the even if the initial data for the PDE (\ref{3.3}) is smooth the solution may cease to be valid later on. In such cases discontinuities appear in the state of the physical system governed by the conservation law (\ref{3.2}), see section 5 in CHAPTER 3 in \cite{pp-rr-book}. Then one has to look for a weak solution of the conservation law (\ref{3.2}). For the definition of a weak solution solution, we refer to any standard text, here we just mention that a weak solution may not be smooth but may have discontinuities almost everywhere and satisfies the integral form (\ref{3.1}).

But it turns out that a weak solution may not be unique due to presence of physically unrealistic discontinuities.  This happened for the discontinuities in gasdynamics in the second half of the 19th century and many important mathematicians were baffled. It was only in 1910 that Lord Rayleigh pointed out that an acceptable discontinuity is the one for which entropy of the fluid elements increase (see the classical book \cite{cou-fri}, Section 51 for a historical account) after crossing the discontinuity. 

It was in 1953 and 1957, Peter lax in USA and in 1957 O. Oleinik in USSR developed the mathematical theory for the conservation law (\ref{3.2}) and defined the generalised solution (now called weak solutions), which include integrable functions with many discontinuities. I give  reference only to Lax (1957) \cite{lax-57}. For a conservation law in which $\frac{\dot F(u)} {\dot H(u)}= a(u)$, say, is a monotonically increasing function of $u$, the Lax condition for an acceptable discontinuity is 
\begin{equation}
\label{3.4}
a(u_r)<S= \dot X(t)<a(u_l),
\end{equation}
where $u_l$ and $u_r$  are the states on the left and right of the discontinuity at $x=X(t)$ respectively and $S=\dot X(t)$ is the velocity of the discontinuity. The jump discontinuity relation is given in terms of the jumps $[H]$ and $[F]$ respectively in $H(u)$ and $F(u)$ by
\begin{equation}
\dot X(t) (H(u_\ell (t)) - H(u_r (t))) \ = \ F(u_\ell) - F(u_r), ~~\textit{or symbolically}~~ \dot X(t) [H]\ = \ [F].
\label{3.5}
\end{equation}

\begin{defn}
The condition (\ref{3.4}) is called the Lax entropy condition.
\end{defn}
\begin{defn}
Stable or entropy satisfying discontinuity is called a \textbf{shock}.
\end{defn}

We now state an important theorem in a very simple language: 

\vspace{-5mm}
\begin{theo} Solution of an initial value problem of the conservation law (\ref{3.2}) with shocks exits for all $t>0$ and is unique.
\end{theo}

\section{System of conservation laws in one-space dimension} KCL is a system of conservation laws. Hence, we consider 
 \begin{equation}
 \label{4.1}
{{\uu H_t(\uu{u})}}+ \uu{F}_{x}(\uu{u}) = 0,
 \end{equation}
where $\uu u = (u_1,u_2,\cdots,u_n), ~\uu H = (H_1,H_2,\cdots,H_n)^T$ and $\uu F = (F_1,F_2,\cdots,F_n)^T$, are vectors with $n$ components.
As in the case of a single conservation law (\ref{3.3}), by (\ref{4.1}) we mean an integral formulation like (\ref{3.1}) and we can define a weak solution (with discontinuities) of (\ref{4.1}) 

Differential form of (\ref{4.1}) is a system of $n$ first order quasilinear PDEs 
 \begin{equation}
 A( x, t, \uu{ u})\uu{ u}_t + B( x, t, \uu{ u})\uu{u}_x = 0
\label{4.2}, 
 \end{equation}
where the $n \times n$ matrices $A$ and $B$ are respectively Jacobian matrix of functions of $\uu H(u)$ and $\uu F(u)$, given by
\begin{equation}
A(\uu{u}) =  \uu \nabla_{\uu{u}}  \uu{H} \equiv [\uu H_{u_1}, \uu H_{u_2}, \cdots, \uu H_{u_n}]
\: \mbox{and} \: B^{(\alpha)}(\uu{ u})= \uu \nabla_{\uu{u}} \uu{F}^{(\alpha)} \equiv [\uu{F}^{(\alpha)}_{u_1}, \uu{F}^{(\alpha)}_{u_2}, \cdots, \uu{F}^{(\alpha)}_{u_n}]. 
 \label{4.3}
 \end{equation} 
We assume that the matrix $A$ is non-singular.

\noindent\textbf{\large{Physical interpretation:}}  One of the most important examples of a hyperbolic system of partial differential
equations is the system governing the gasdynamic equations. The theory of shock wave arose from the work of many important mathematicians (\cite{cou-fri}, Section 51). Some of them are: Stokes (1848), Riemann (1860), Rankine and Hugoniot (1869), Hugoniot (1887), Rayleigh (1910) culminating in the work of Lax (1957). We consider here of Euler's equations governing the motion of a polytropic gas for
\begin{equation}
  \label{4.8}
  {\bf u} = (\rho, q, p)   
\end{equation}
where $\rho $ is the mass density, $q$ the fluid velocity, $p$ the pressure
and $\gamma$ is the ratio of specific heats assumed to be constant.
These equations are
\begin{align}
  \rho_t + q \rho_x + \rho q_x &= 0,  \label{4.9}\\
  q_t + qq_x + \frac{1}{\rho} p_x &= 0, \label{4.10} \\
  p_t + q p_x + \gamma p q_x &= 0  \label{4.11}.
\end{align}
The eigenvalues are
\begin{equation}
  \label{4.12}
  c_1 = q - a, \quad c_2 = q, \quad c_3 = q + a,
\end{equation}
where $a$ is the local sound velocity in the fluid relative to the fluid
particles and is given by
\begin{equation}
  a^2 = \gamma p / \rho.  \label{4.13}
\end{equation}
When a shock appears in the third characteristic field $c_3$, it is a forward facing shock and the inequality  (\ref{3.4}) implies 
\begin{equation}
q(\uu u_r) + a(\uu u_r) < S< q(\uu u_\ell) + a(\uu u_\ell).
\end{equation}
This result has a very important \textbf{physical interpretation:} A shock moves with a supersonic speed with respect to the state ahead of it and with a subsonic speed with respect to the state behind. It also implies that the specific entropy of the fluid elements increase after it passes through a shock. 

\section{Ray coordinates $(\xi,t)$ for an isotropic motion of a moving curve $\Omega_t$ in 
$(x,y)$-plane:} Associated with a moving curve there exist rays, which carry different points $(x_0,y_0)$ of the $\Omega_0$ as the curve moves. For a curve evolving isotropically, the ray velocity $\uu \chi = (\chi_1, \chi_2)$ is in the normal direction $\uu n=(n_1,n_2)$ of $\Omega_t$. Let the surface $\Omega_t$ be represented by 
\begin{equation}
\label{5.1}
\Omega_t: \varphi(x,t) =0
\end{equation}
 and $m$ be a suitably normalised velocity of $\Omega_t$, then 
\begin{equation}
\label{5.2}
\uu \chi = (\chi_1, \chi_2), ~~\textit{where}~ \chi_1=n_1m, \chi_2=n_2m.
\end{equation}
Note that the $m$ and $\uu n$ are given by
\begin{equation}
\label{5.3}
C \equiv m =-\frac{\phi_t}{\mid \bigtriangledown \varphi \mid}, ~~~~~~~~~
{\uu{n} =\bigtriangledown \varphi / \mid \bigtriangledown \varphi \mid}
\end{equation} 
We write a parametric equation of $\Omega_{t}: \varphi(x,y,t)=0$ in the form
\begin{equation}
\Omega_{t} : x=x(\xi, t), \:\:  y=y(\xi,t),
\label{5.4}
\end{equation}
where $\xi$ the parametric variable on the curve $\Omega_t$ at time $t$. From (\ref{5.3}) we get the eikonal equation
\begin{eqnarray}
\label{5.5}
\varphi_{t} + m\mid\bigtriangledown \varphi \mid =0
\end{eqnarray}
\noindent which is a first order nonlinear partial differential
equation for $\varphi$ describing evolution of $\Omega_{t}$. 

\noindent\textbf{Ray coordinates:} $(\xi,t)$ forms ray coordinates in $(x,y)$-plane such that in (\ref{5.4}) constant values of $t$ give the positions of the propagating curve $\Omega_{t}$ at different times and $\xi =$ constant represents a ray. The coordinate $\xi$ is not unique and can be replaced by any monotonic function of it. Consider now a domain $\mathcal{D}_p$, swept by $\Omega_t$, of the $(x,y)$-plane. Then there is one to one mapping between $\mathcal{D}_p$ in $(x,y)$-plane and a domain $\mathcal{D}_r$ in $(\xi,t)$-plane\footnote{The subscript $p$ refers to the physical plane and $r$ to the ray coordinate plane}. 
\begin{figure}
    \centering
    \includegraphics[height=8cm]{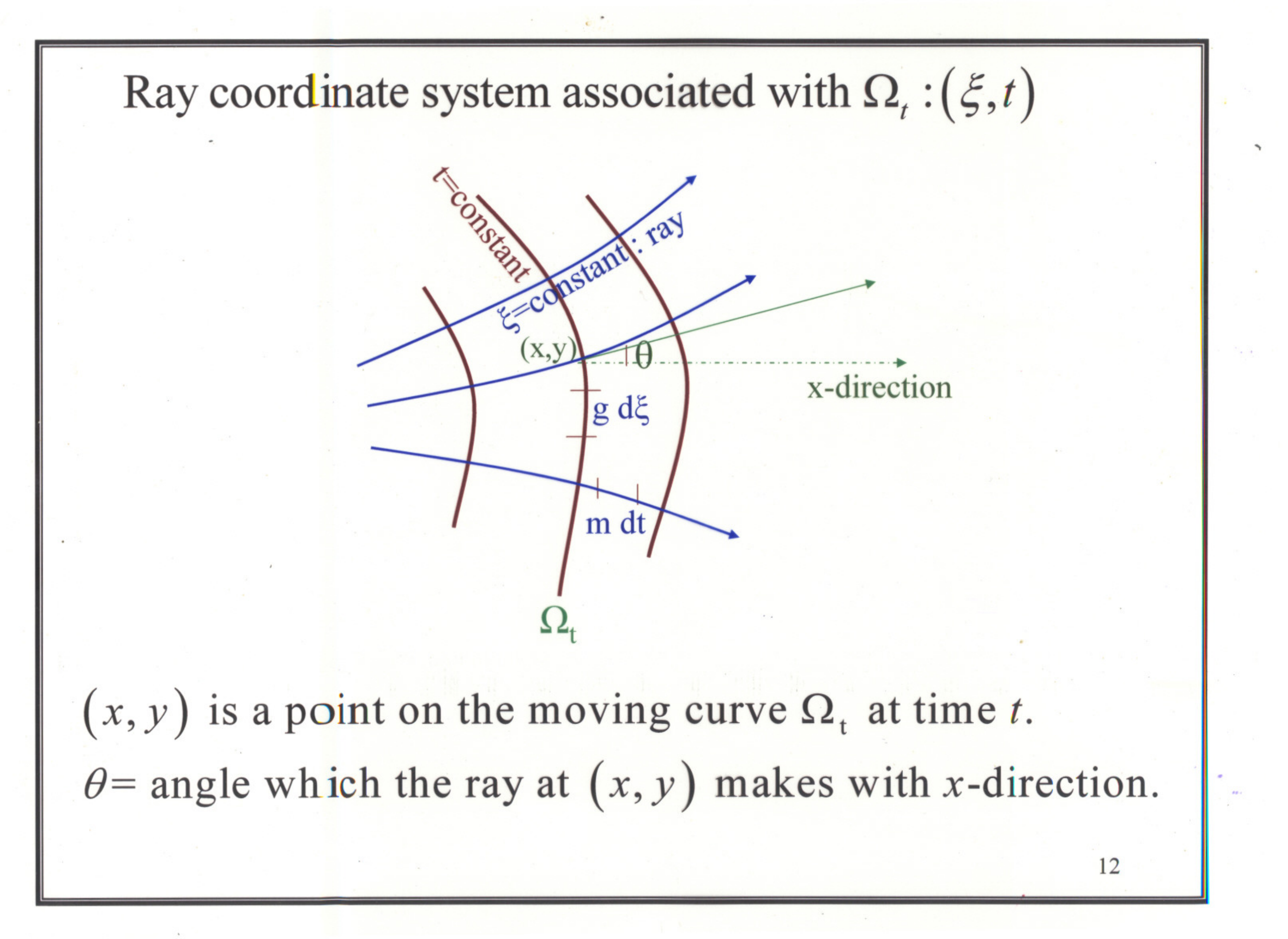}
\caption{\footnotesize In $(x,y)$-plane, $\xi$ = constant is a ray and $t$=constant is $\Omega_t$.}
  \end{figure}
Let $\theta$ be the angle which the ray at a point $(x,y)$ makes with the $x$-axis. 
The time rate of change along a ray is 
\begin{equation}
\label{5.6}
\frac{d}{dt} = \frac{\partial}{\partial t} + m\left(\cos \theta \frac{\partial}{\partial x}+\sin \theta \frac{\partial}{\partial y}\right).
\end{equation}
The tangential derivative $\frac{\partial}{\partial \xi}$ along the curve $\Omega_t$ is given by
\begin{equation}
\label{5.7}
\frac{\partial}{g\partial \xi} = -\sin \theta \frac{\partial}{\partial x} + \cos \theta \frac{\partial}{\partial y },
\end{equation} 
where $g$ is the metric associated with the coordinate $\xi$, i.e $gd\xi$ is an element of the distance along $\Omega_t$. Note that $m$ is the velocity of a point on the ray, i.e.\ $mdt$ is an element of the distance along a ray, as shown in the Figure 2.
In ray coordinates the time rate of change along a ray and rate of change along  $\Omega_t$ are denoted by
\begin{equation}
\label{5.8}
\frac{d}{dt} =\partial_t ~~~~~~~~~~~\textit{and} ~~~~~~~~~~~ \frac{\partial}{\partial \xi}=\partial_\xi.
\end{equation}
We can use the theory of a first order PDE, where the characteristic ODEs of (\ref{5.5}) give:
\begin{equation}
\label{5.9}
\frac{dx}{dt}= m \cos \theta, ~~~ \frac{dy}{dt}= m \sin \theta
\end{equation}
and
\begin{equation}
\label{5.10}
\frac{d\varphi_{x}}{dt}= -\mid\bigtriangledown \varphi \mid m_{x},~~~~ \frac{d\varphi_{y}}{dt}= -\mid\bigtriangledown \varphi \mid m_{y}.
\end{equation}
We shall derive an equation for $\theta$ from any one of the two equations in (\ref{5.10}).

Differentiating $(n_1 \equiv \cos \theta )= \frac{\varphi_x}{\mid\bigtriangledown \varphi \mid}$ along a ray and using $1-n_1^2 = n_2^2$ we get
\begin{equation}
\label{5.11}
\frac{dn_1}{dt} = \frac{n_2}{\mid\bigtriangledown \varphi \mid}\left\lbrace n_2 \frac{d\varphi_{x}}{dt} -n_1 \frac{d\varphi_{y}}{dt} \right\rbrace.
\end{equation}
Using (\ref{5.10}) we finally get
\begin{equation}
\label{5.12}
\frac{d \theta}{dt} = 
\left(-\sin \theta \frac{\partial}{\partial x} + \cos \theta \frac{\partial}{\partial y}\right)m.
\end{equation}
Equations (\ref{5.9}) in ray coordinates become
\begin{equation}
\label{5.13}
x_t= m \cos \theta, ~~~ y_t= m \sin \theta.
\end{equation}
Once the parametric representation of (\ref{5.4}) of $\Omega_t$ is made,
the metric $g$ is given by
\begin{equation}
\label{5.14}
 g^2 = x_\xi^2 + y_\xi^2.
\end{equation}
Differentiating partially with respect to $t$ and using (\ref{5.13}) and
\begin{equation}
\label{5.15}
\frac{x_\xi}{g} = -\sin \theta, ~~~~~~ \frac{y_\xi}{g} = \cos \theta,
\end{equation}
we can derive 
\begin{equation}
\label{5.16}
 g_t = m \theta_\xi.
\end{equation}
Equations (\ref{5.12}) in ray coordinates become 
\begin{equation}
\label{5.17}
\theta_t = - \frac{1}{g} m_\xi
\end{equation} 

The results contained in equations (\ref{5.16}) and (\ref{5.17}) are physically realistic. We verify for an expanding circular curve $\Omega_t$ with centre at the origin with a constant velocity $m>0$ originally. Then it remains circular all the time. Since $m_\xi=0$, (\ref{5.17}) implies that the $\theta$ remains constant along a ray, which is a straight line. Since $\theta_\xi>0$, (\ref{5.16}) shows that $g$ increases with $t$. This also consistent with the fact that part 
$\int_{\xi_\ell}^{\xi_r}gd\xi$ of $\Omega_t$ keeps on increasing in the fixed interval
$({\xi_\ell},{\xi_r})$. Any formulation of KCL should lead to (\ref{5.16}) and (\ref{5.17}) for a smooth evolution of the moving curve.

\section{Two-dimensional ($2$-D) system of kinematical conservation laws (KCL):} We shall give formulation of $2$-D KCL in terms of the density of ``KCL conserved variable" and ``KCL flux". Kinematics refers to study of the motion of bodies without reference to mass or force. KCL involves purely geometric objects and is in $(\xi,t)$-plane.

\noindent\textbf{A proposal for density of KCL conserved variable:} Associated with $\Omega_t$, there are two important vectors \\
(1) ray velocity $(m \cos \theta, m\sin\theta)$ (see (\ref{5.13})) normal to the curve $\Omega_t$ and \\
(ii) a velocity $(-g \sin \theta, g\cos\theta)$ (see (\ref{5.15})) tangential to the curve $\Omega_t$. \\
We propose that the conserved variable in KCL is the tangential velocity mentioned in (2) and the flux vector is the ray velocity mentioned in (1). 

The {balance between} time rate of change of total  conserved quantity (now a tangential vector with two components on $\Omega_t$) from the point identified by $\xi_\ell$ on $\Omega_t$ to that identified by $\xi_r$ on $\Omega_t$ and the flux from the two ends is expressed as a conservation law: 
\begin{equation}
\label{6.1}
\frac{d}{dt} \int_{\xi_ell}^{\xi_r}(-g\sin \theta, g\cos \theta)(\xi,t)d\xi  =  \left\lbrace(m\cos \theta, m\sin \theta)\right\rbrace(\xi_\ell,t)-
\left\lbrace(m\cos \theta, m\sin \theta)\right\rbrace(\xi_r,t).
\end{equation}

If $\theta$ is constant on $\Omega_t$, then the front is a straight line at time $t$. If $m$ is constant on it, then every point at time $t$ moves with the same velocity. When both are constant on $\Omega_t$, the curve propagates as a straight line parallel to itself and the equation (\ref{6.1}) shows that the $(-g\sin \theta, g\cos \theta)(\xi,t)$ on $\Omega_t$ between the points corresponding $\xi_\ell$ and $\xi_r$ remains constant as time evolves, i.e.\ \textbf{it remains conserved}. 

Symbolically the pair conservation laws (\ref{6.1}) of $2$-D KCL is denoted (see footnote (4) related to ({\ref{3.2}))  by
\begin{align}
(g\sin\theta)_t + (m\cos\theta)_\xi &= 0, \label{6.2} \ \mbox{and} \\
(g\cos\theta)_t - (m\sin\theta)_\xi &= 0. \label{6.3}
\end{align}
We wrote the form of KCL in our 1992 article (by Morton, Prasad and Ravindran) without mentioning the density vector of the conserved quantity and the flux vector.
The differential form of KCL (\ref{6.2}) and ((\ref{6.3}) are (\ref{5.16}) and (\ref{5.17}), which we derived from the ray equations. Let us notice the according to the comment made after (\ref{5.17}), the $2$-KCL is a physically realistic system.


\section{Theorem: KCL and ray equations are equivalent for the evolution of a smooth $\Omega_t$.}
\textbf{KCL implies the ray equation:}
First note that $\theta$ in KCL is the same as that in the ray equations (\ref{5.9}) and (\ref{5.12}).
According to the fundamental integrability  theorem (\cite{cou-joh}-page 104), the conditions (\ref{6.2}) and (\ref{6.3}) imply the existence of a vector $\uu{\tilde x}$ satisfying
\begin{equation}
\label{6.7}
\uu{\tilde x}_t= m(\cos \theta,\sin \theta) ~~\textit{and}~~  \uu{\tilde x}_\xi=g(-\sin\theta,\cos \theta).
\end{equation}
The normal velocity of $\Omega_t$ is $m$ which must be equal to that in the ray equations 
$-\varphi_t / |\uu \nabla_{(\tilde x, \tilde y)} \varphi|$. The function $\varphi$ appearing in (\ref{5.1}) satisfies the eikonal equation (\ref{5.5}) with $x$ and $y$ replaced by $\tilde x$ and $\tilde y$. The characteristic equations of this eikonal equation give the equations (\ref{5.9}) and (\ref{5.12}) again with $x$ and $y$ replaced by $\tilde x$ and $\tilde y$. Now we need not distinguish between $\uu {\tilde{x}}$ and  $\uu x$. 

Thus one part of the equivalent is established.

\noindent\textbf{The ray equations imply KCL:} We have identified a domain $\mathcal{D}_p$ in $(x,y)$-plane and its image domain $\mathcal{D}_r$ in $(\xi,t)$-plane, see the paragraph just above (\ref{5.6}). Increments $d\xi$ and $dt$ in the ray coordinates leads to a displacement $dx$ and $dy$ in $(x,y)$-plane given by (see Figure 2)
\begin{equation}
\label{7.2}
dx=\cos \theta (mdt)- \sin \theta (gd\xi),~~~~~dy = sin \theta(mdt) + cos \theta (gd\xi).
\end{equation}
These two results also follow by taking projections of line elements $mdt$ and $gd\xi$ on the $x$ and $y$ axes.

We use these two to  equate $x_{\xi t} = x_{t\xi}$ and $y_{\xi t} = y_{t\xi}$ and derive the pair of relations (\ref{6.2}) and (\ref{6.3}) in the conservation form. 
Thus another part of the equivalent is also established. 

\section{Three-dimensional ($3$-D) system of kinematical conservation laws (KCL):}
Unlike the $2$-D KCL the theory of $3$-D KCL is quite involved. We describe the results briefly, the details of its theory and applications are available in papers \cite{aru-pp-09}, \cite{aru-jor-mod}, \cite{aru-pp-12} and \cite{pp-16b}. All results on on ray theory, hyperbolic conservation laws and KCl are available in \cite{pp-book-18}, first 12 pages are available at $http://www.math.iisc.ernet.in/\sim prasad/prasad/book/books.html$..

\noindent\textbf{Ray equations in $3$-D space.} The results (\ref{5.3}) are true in space of any dimensions, hence (\ref{5.5}) is the eikonal equation for $\Omega_t: \varphi(\uu x,t)=0$ in $3$-D is also. We denote the components of the unit normal to $\Omega_t$ as
\begin{equation}
\label{8.1}
\uu n = (n_1,n_2,n_3),\: \: \left|\uu{n}\right|=1.
\end{equation}
In this case characteristic equations (as in the case (\ref{5.9}) and (\ref{5.12})) ray equations  take simple forms
\begin{align}
\frac{d\uu{x}}{d t}&=m\uu{n}, \label{8.2}\\
\frac{d \uu{n}}{d t} &=-\uu{L}m:=-\left(\uu \nabla-\uu{n}
  \langle\uu{n},\uu \nabla\rangle\right)m. \label{8.3}
\end{align}
The operator $\uu L$ defined above is obtained by subtracting from the gradient $\uu \nabla$, its component in normal direction $\uu n$ and hence it represents a tangential derivative on $\Omega_t$. 

\noindent \textbf{Ray coordinate system:}\index{ray coordinates}
Consider the isotropic evolution\footnote{In an isotropic evolution each point of the surface moves in the normal direction of $\Omega_t$.} of a $2$-D surface $\Omega_t$ in $\RR^3$. Let $\uu{\xi}=(\xi_1,\xi_2)$ be a set of surface coordinates, which also evolve with time $t$ (see Figure 3).
 \begin{figure}
    \begin{center}
      \includegraphics[height=8cm]{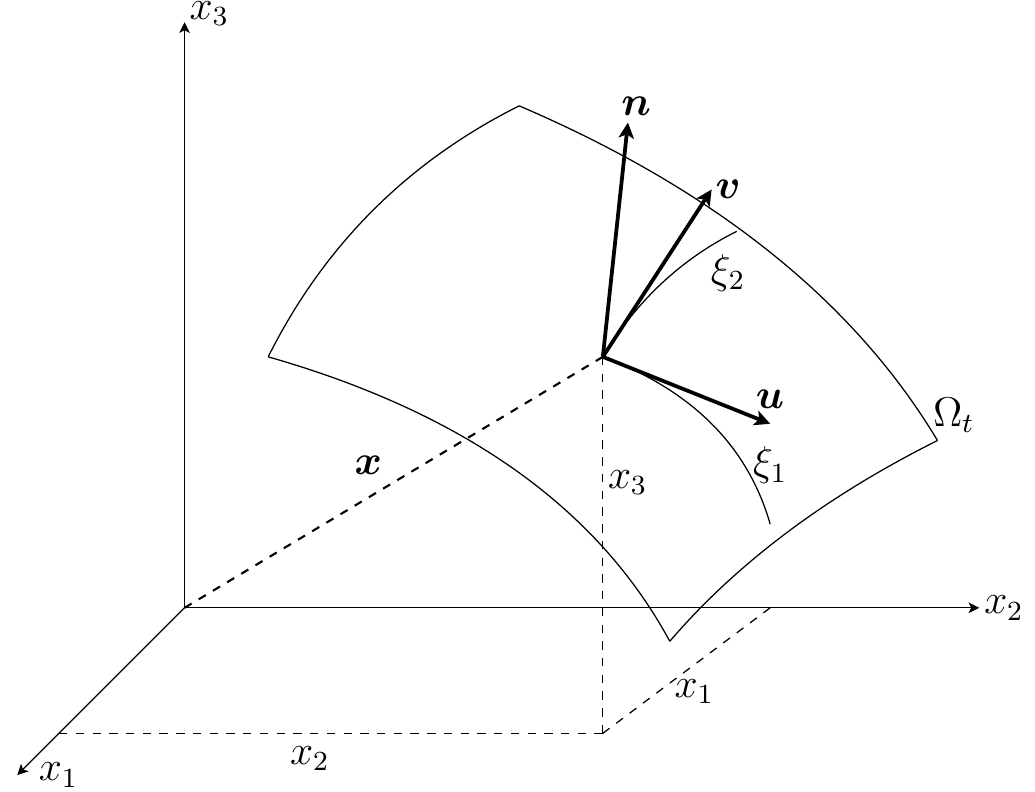}
    \end{center}
    \caption{\footnotesize A ray coordinate system of a surface $\Omega_t$.} 
  \end{figure}
At a fixed time $t$ the surface $\Omega_t$ in $\uu{x}$-space is generated by a $2$ parameter family of curves, the parameters being $\xi_1$ and $\xi_2$. Along a member of the first family of curves coordinate $\xi_1$ varies and the parameters
$\xi_2$ is constant. Similarly along a member of the second family of curves coordinate $\xi_2$ varies and the parameters $\xi_1$ is constant.  Through each point $\uu{\xi}$ on $\Omega_t$ there passes a ray in normal direction $\uu n$. Let 

Let $\uu u$ and $\uu v$ be unit vectors along the $\xi_1$ and $\xi_2$ families of the coordinates, then 
\begin{equation}
\label{8.4}
\uu{n}=\frac{\uu{u}\times\uu{v}}{\left|\uu{u}\times\uu{v}\right|}.
\end{equation}
\noindent\textbf{Density of the conserved vector and flux vectors for $3$-D KCL:} Let the metric associated with the surface coordinate $\xi_p$ be $g_p$, $p=1,2$. Then $g_p d\xi_p$ (no sum over the repeated subscript here) is an element of distance along the coordinate line along which $\xi_p$ varies. The speed of a point moving along the ray with the ray velocity is $m$ then, while moving with the the ray velocity, $mdt$ is the displacement along the ray in time $dt$. Thus $m$ is the metric associated with the coordinate $t$. 

As in the case of the $2$-D KCL we take the density of the conserved vector  along the first family of coordinates to be $g_1 \uu u $  and that along the second family of coordinates to be $g_2 \uu v $. We also propose that the flux vectors along both families to be the ray velocity $m\uu n$. 

Similar to the case of $2$-D ray coordinates as mention in the paragraph above the equation (\ref{5.6}), we note that there is one to one mapping between $\mathcal{D}_p$ in $(x,y,z)$-space and a domain $\mathcal{D}_r$ in $(\xi_1,\xi_2,t)$-space.

We first consider the formulation in the subspace of $(x_1,x_2,x_3)$-space in which $\xi_2$ is constant. This means we consider the formulation in $(\xi_1,t)$-plane of $(\xi_1,\xi_2,t)$-space. The conserve density vector is $g_1\uu{u}$ and the flux is $m\uu n$.  As in the case of $2$-D KCL, the integral formulation of the first set of conservation law is 
\begin{equation}
  \label{8.5}
\frac{d}{dt} \int_{{\xi_1}_\ell}^{{\xi_1}_r} g_1(\xi_1,\xi_2,t) \uu u_1(\xi_1,\xi_2,t)d\xi_1  =m({\xi_1}_\ell,\xi_2,t) \uu n({\xi_1}_\ell,\xi_2,t) - m({\xi_1}_r,\xi_2,t)) \uu n({\xi_1}_r,\xi_2,t)), \ \xi_2=\mbox{const.}
\end{equation}
Then we consider the formulation in $(\xi_2,t)$-plane of $(\xi_1,\xi_2,t)$-space. The conserve density vector is $g_2\uu{u}$ and the flux is $m\uu n$.  Now we write the second set of conservation law as
\begin{equation}
  \label{8.6}
\frac{d}{dt} \int_{{\xi_2}_\ell}^{{\xi_2}_r}g_2(\xi_1,\xi_2,t)\uu u_2(\xi_1,\xi_2,t)d\xi_2  =(m (\xi_1,{\xi_2}_\ell,t)\uu n(\xi_1,{\xi_2}_\ell,t) -(m(\xi_1,{\xi_2}_r,t)\uu n(\xi_1,{\xi_2}_r,t)), \ \xi_1=\mbox{const.}
\end{equation}
The symbolic form of the conservation laws (\ref{8.5}) and (\ref{8.6}) are 
\begin{equation}
(g_1\uu{u})_t-(m\uu{n})_{\xi_1}=0, ~~~~
(g_2\uu{v})_t-(m\uu{n})_{\xi_2}=0 \label{8.7}.
\end{equation}
$\uu u$ and $\uu v$ each has three components and hence the $3$-D KCL has 6 equations. But since these two vectors are unit vectors only $4$ of these are independent. 

\noindent\textbf{\begin{theo} $3$-D KCL and ray equations (\ref{8.2}) and (\ref{8.3}) are equivalent for the evolution of a smooth $\Omega_t$.\end{theo}} 
We have proved this result in Section 7 for $2$-KCL, in we omit the proof here, since different proofs are available in references quoted in the beginning of this section.

An explicit differential form of the the $3$-D KCL are available in Section 4 of the paper \cite{aru-pp-09}, a preprint of which is available at \\
($http://math.iisc.ernet.in/\sim prasad/prasad/preprints/3d\_kcl.pdf$). 

\noindent\textbf{Geometrical solenoidal constraint:} This is an interesting aspect of KCL in three and higher dimensions. We derive it here for $3$-D KCL. The
displacement $d \uu x$ in $\uu x$-space due to increments
$d\xi_1$, $d\xi_2$ and $d t$ in $(\uu \xi,t)$-space is given by 
\begin{equation}
\label{8.8}
d \uu x=(g_1\uu u)d\xi_1+(g_2\uu v)d\xi_2+(m\uu n)d t.
\end{equation} 
For a smooth moving surface
$\Omega_t$, we equate $\uu x_{\xi_1t}=\uu x_{t\xi_1}$ and
$\uu x_{\xi_2t}=\uu x_{t\xi_2}$, and get the 3-D KCL (\ref{8.7}).
We also equate $\uu x_{\xi_1\xi_2}=\uu x_{\xi_2\xi_1}$ and
derive 3 more scalar equations contained in
\begin{equation}
\label{8.9}
(g_2\uu v)_{\xi_1}-(g_1\uu u)_{\xi_2}=0.
\end{equation}
(\ref{8.9}) is purely geometrical results as the $m$ does not appear in it and we call it \textbf{geometrical solenoidal constraint}.
From the equations (\ref{8.7}) we can show that
$(g_2\uu v)_{\xi_1}-(g_1\uu u)_{\xi_2}$ does not depend on
$t$. If any choice of coordinates $\xi_1$ and $\xi_2$ on
$\Omega_0$ the condition (\ref{8.9}) is satisfied at
$t=0$ then it follows that (\ref{8.9}) is automatically satisfied on $\Omega_t$ for all $t>0$. 

\noindent\textbf{Control of Jordan mode by geometrical solenoidal constraint:} We have not yet mentioned some applications of the KCL, which have been done by our group. If we consider KCL based weakly nonlinear ray theory (WNRT) or shock ray theory (SRT) ( see \cite{pp-book}, \cite{bas-pp-05} and any of the references mentioned above) the $2$-D problems remain quite simple but $3$-D problems become quite complex. For example the $3$-D WNLRT equations become a degenerate with a multiple eigenvalue for which the eigenspace is incomplete. In such a case Jordan mode generally appears in numerical solution which grow in time. It is not the place to discuss the details for controlling the Jordan mode for which we refer to \cite{aru-jor-mod}. 

\section{Some applications: KCL based weakly nonlinear ray theory (WNLRT) and KCL based shock ray theory (SRT)} The system of KCL equations is under-determined. 
For example the $2$-D KCL equations are two in number for three dependent variables $m,\theta$ and $g$. 
The $3$-D KCL equations are six in number for seven dependent variables (two components of $\uu u$, two components of $\uu v$, $g_1$, $g_2$ and $m$). 
When the moving curves represents a physically observable curve such as a wavefront or a shock front or crest line of a curved solitary wave on the surface of water \cite{bas-pp-03}, we hope to close the KCL system of equations. We do not go into details but just reproduce figures from the references mentioned above.

\subsection{Exact solutions showing resolution of a caustic in $2$-D by WNLRT:} {Figure on the right is reproduced from \cite{pp-book}}.
\begin{figure}[hbt]
\hspace{10.0mm}
\vspace{-10mm}
\begin{minipage}[s]{2.3in}
\includegraphics[height=0.3\textheight]{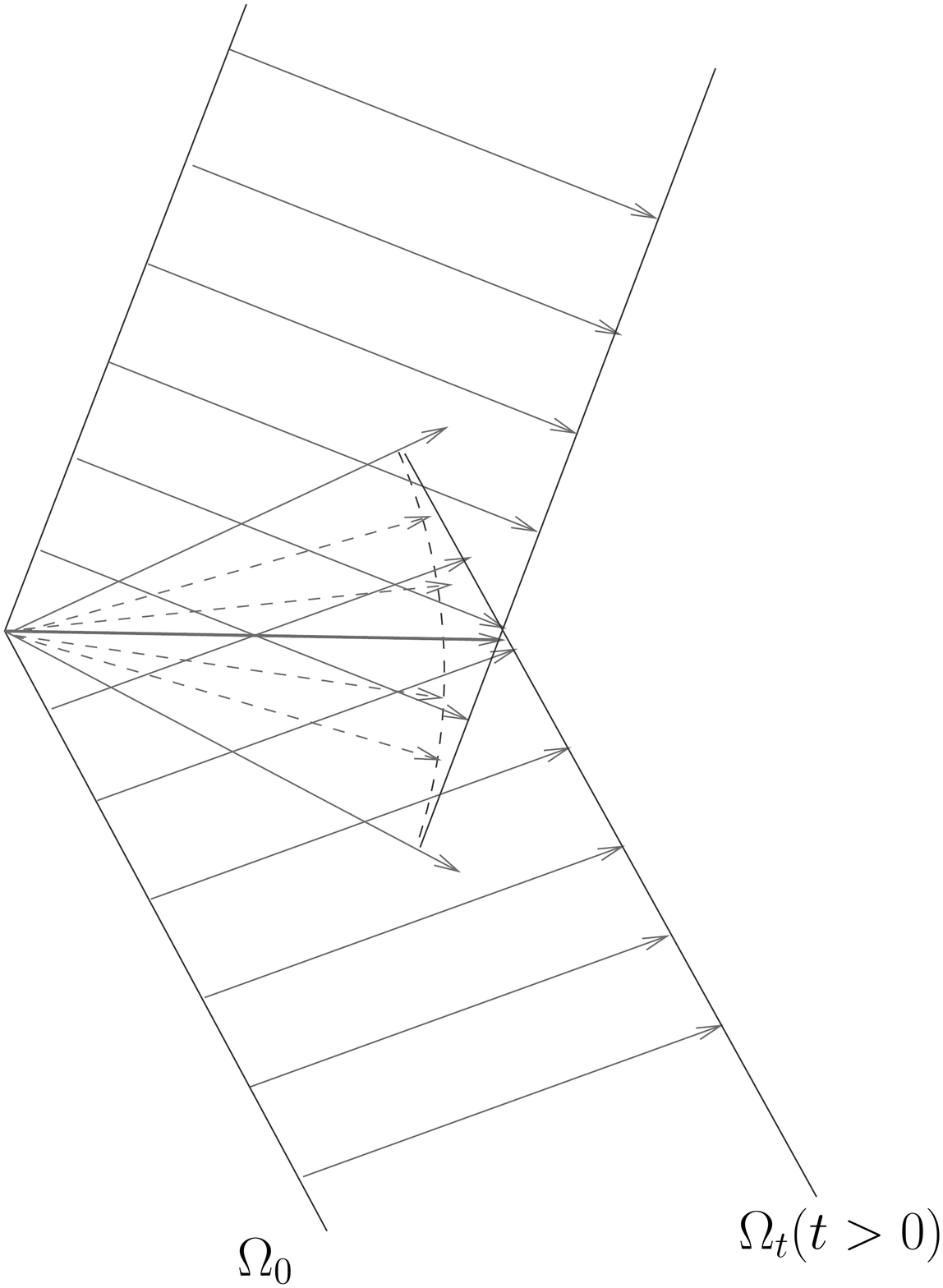}
\end{minipage} 
\begin{minipage}[s]{2.3in}
\includegraphics[height=0.3\textheight]{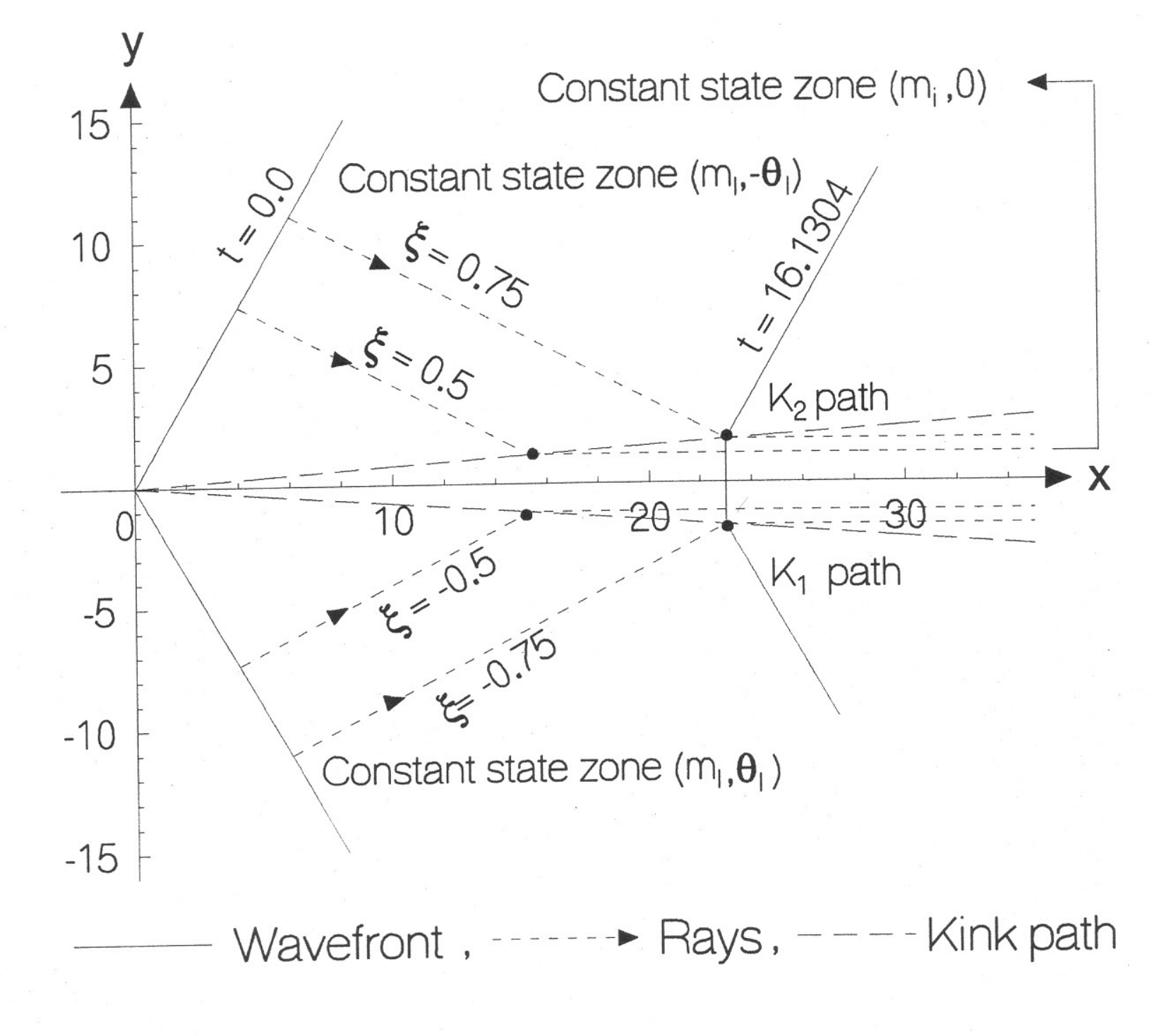}
\end{minipage}
\vspace{10mm}
\caption{{\footnotesize Figure on the left has a caustic: Linear wavefront which starts with a concave wedged shape. \textcolor{blue}{.....} wavefront produced by the corner O  by Huygen's method. This is an arc of a circle with centre at O.
 \textcolor{blue}{-----} wavefront produced by the smooth part of  the initial wavefront by ray theory or Huygen's method. ~~~~~~~~~~~~~~~~ {In the Figure on the right, the caustic of the linear theory is completely resolved by WNLRT. This is an exact solution of WNLRT equations. We find that two kinks (a new type of singularities shown by dots) have appeared. Across a kink directions of the tangents of the wavefront and directions of the rays suffer jumps. These features can be captured only by a KCL based theory.}}}
\end{figure} 

\newpage
\subsection{Evolution of a shock front of periodic shape in two space dimensions:}
{The picture is reproduced from a paper by A. Monica and PP, 2001, Journal of Fluid Mechanics.}

For a shock front, we denote normal velocity of the shock by $M$ (this means we replace $m$ by $M$). We take the initial shock front $\Omega_0$ to be in a periodic sinusoidal
shape
\begin{equation}
 x=0.2-0.2 \cos (\frac{\pi y}{2}).
\label{9.1}
\end{equation}
We choose the initial shock Mack number (appropriately non-dimensionalised  shock speed) to be uniform and equal to $M_0 = 1.2$. The shape of the shock at various times are given in the Figure 5. The dots represent kinks which appear later. Such a long time tracking of a shock front with kinks is not easy by a method other than a KCL based formulation. 
\begin{figure}
    \centering
    \includegraphics[height=10cm]{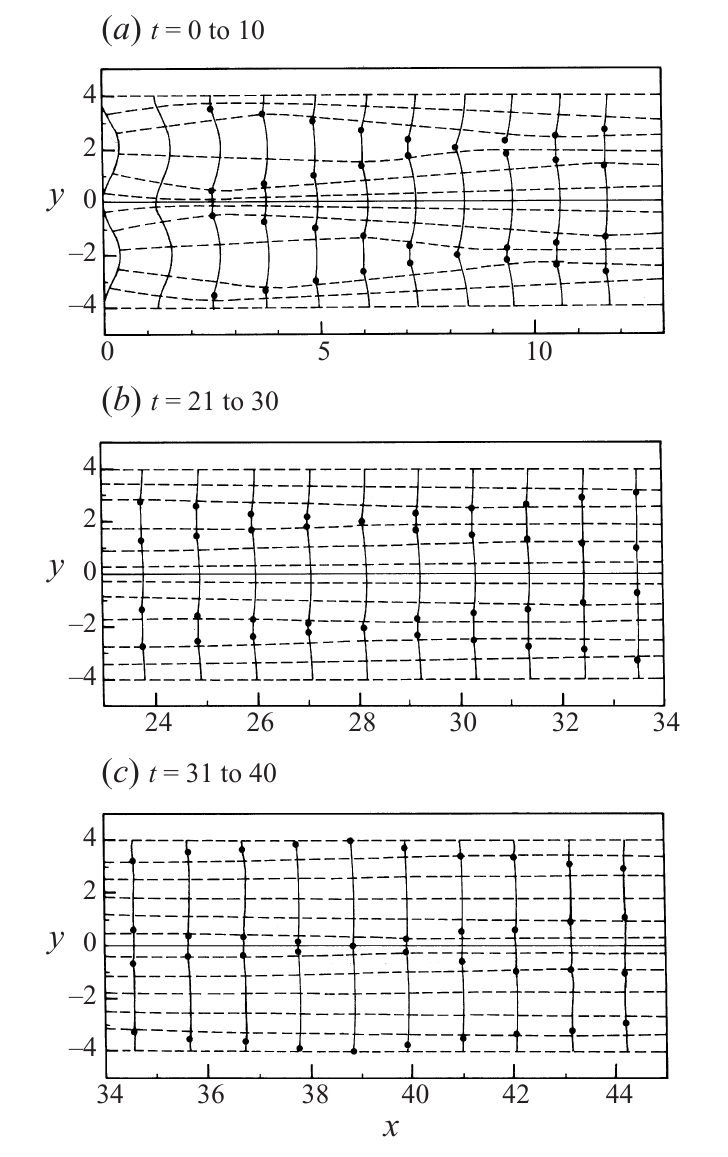}
\vspace{-7mm}
\caption{{\footnotesize Successive positions of an initially sinusoidal shock
front (shown by continuous line) plotted at $t$ = 0, 1, 2, 3, $\cdots$, 40 and rays (shown by broken lines). Kinks appear after $t=1$ and are shown by dots. In the figure we see a kink, say $\mathcal{K}_3$, from a lower period moving upwards and another kink, say $\mathcal{K}_1$, from upper period moving downwards resulting in interaction results $\mathcal{K}_3 \mathcal{K}_1 \rightarrow \mathcal{K}_1 \mathcal{K}_3$ happening many times. The shock has become almost straight and rays parallel to $x$-axis from $t=31$ to $40$. This figure have been reproduced from \cite{mon-pp-01}.}}
\end{figure}

\newpage
\subsection{Evolution of a shock front of periodic shape in three space dimensions:}
{Reproduced from \cite{aru-pp-12}.}
 
 We consider here the initial shock front $\Omega_0$ to be of a periodic shape in $x_1$- and
$x_2$-directions
\begin{equation}
  \Omega_0\colon ~~x_3=\kappa\left(2
    -\cos\left(\frac{\pi x_1}{a}\right)-\cos\left(\frac{\pi
        x_2}{b}\right)\right)
  \label{9.2}
\end{equation}
with the constants $\kappa=0.1,a=b=2$. In Figure 6 we give the plot of the initial shock
front $\Omega_0$, which is a smooth pulse without any kink lines. We choose the Mach number to uniform equal to $1.2$ to be uniform on $\Omega_0$. Though the initial shock front is smooth, as the time evolves, a number of kink lines appear in each period. 

Now we describe the very interesting process of interaction of these kink lines. We proceed to do it with a number of plots of the shock front at different instances.

In Figure 7 we give the surface plots of the shock
front $\Omega_t$ at times $t=10,20,30,40,50,60$ in two periods in each of $x_1$ and $x_2$-directions. As mentioned above, the initial shock
front is smooth, with no kink lines. The front $\Omega_t$ moves up in
the $x_3$-direction and develops several kink lines. Four kink lines
parallel to $x_1$-axis and four parallel to $x_2$-axis can be seen in
the figures on the shock front at times $t\geq 10$. These kink lines
are formed at a time before $t=10$, say about $t=2$.  

The shock front at $t=60$ has almost become plane.
\begin{figure}
  \centering
  \includegraphics[width=0.4\textwidth]{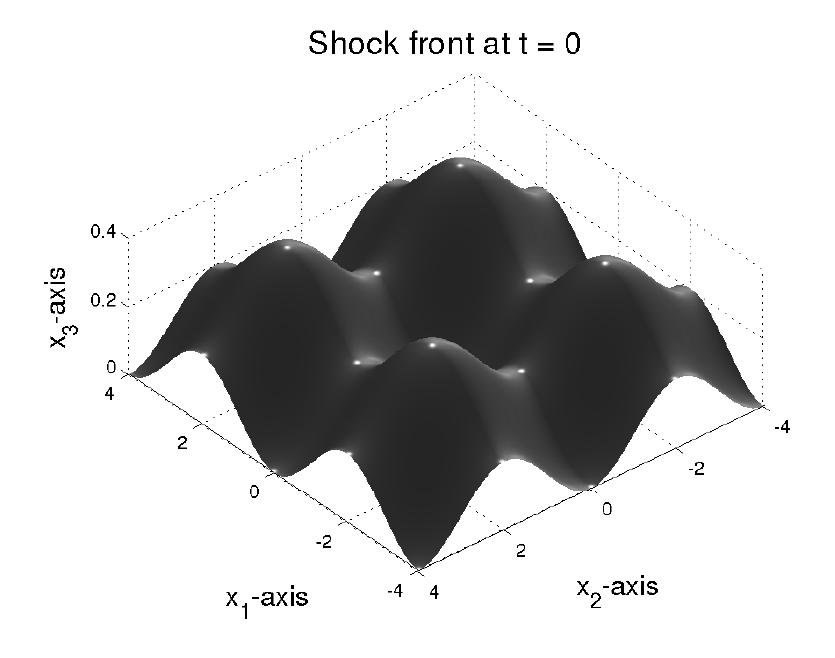}
  \caption{\footnotesize Initial shock front in the shape of a smooth periodic
    pulse.}
  \label{fig10.2}
\end{figure}
\begin{figure}
 \centering
  \includegraphics[width=0.9\textwidth]{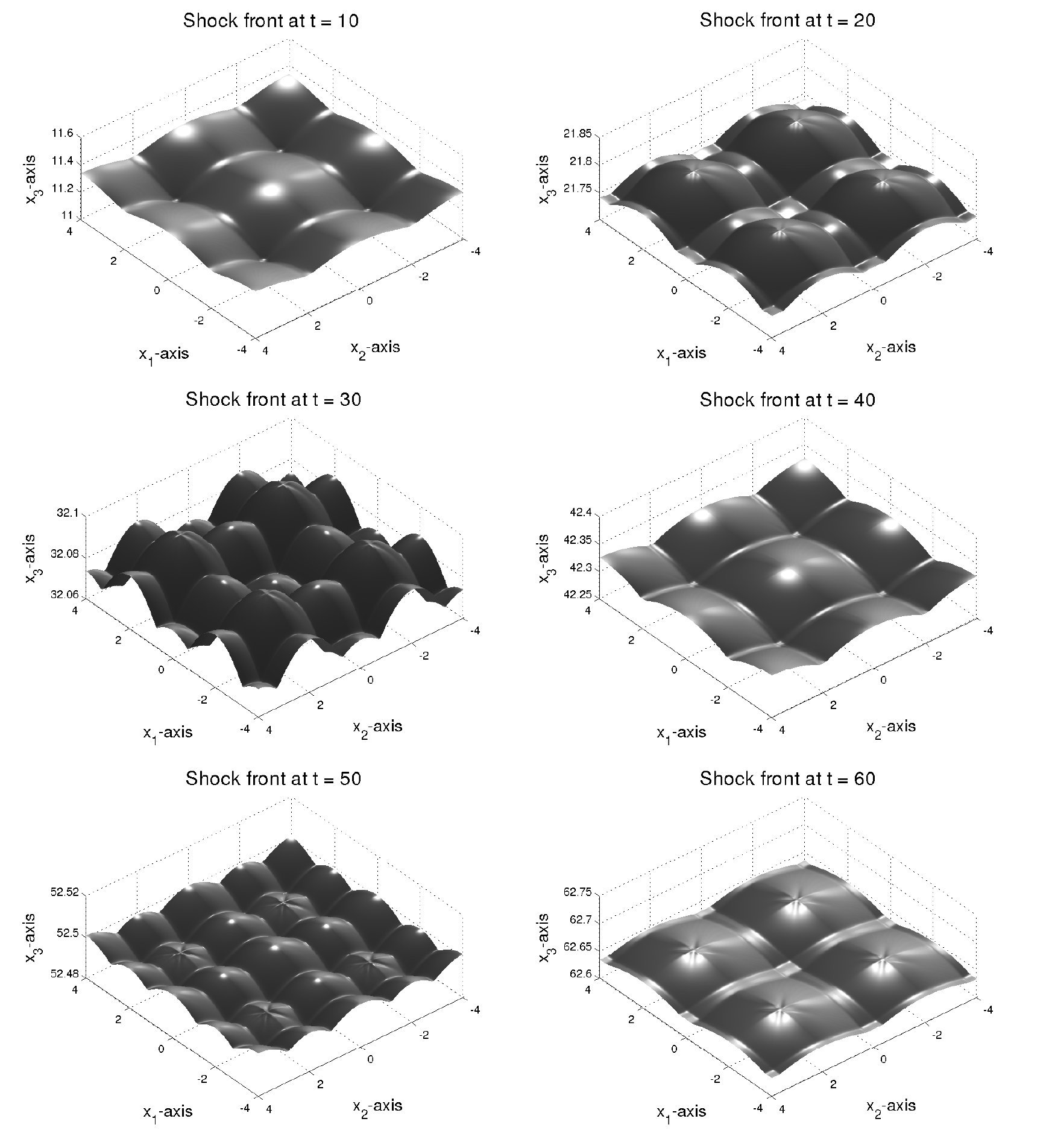}
 \caption{\footnotesize Shock front $\Omega_t$ starting initially in a periodic shape with $M_0=1.2$.  The shock front develops a complex pattern
 of kinks and ultimately becomes planar.}
  \label{fig10.3}
\end{figure}

\newpage 
%

\end{document}

%% file: psfigfx.tex
\def\PsfigVersion{1.9}
\ifx\undefined\psfig\else \fi

\let\LaTeXAtSign=\@
\let\@=\relax
\edef\psfigRestoreAt{\catcode`\@=\number\catcode`@\relax}
\catcode`\@=11\relax
\newwrite\@unused
\def\ps@typeout#1{{\let\protect\string\immediate\write\@unused{#1}}}
\ps@typeout{psfig/tex \PsfigVersion}

\def\figurepath{./}

\def\@nnil{\@nil}
\def\@empty{}
\def\@psdonoop#1\@@#2#3{}
\def\@psdo#1:=#2\do#3{\edef\@psdotmp{#2}\ifx\@psdotmp\@empty \else
    \expandafter\@psdoloop#2,\@nil,\@nil\@@#1{#3}\fi}
\def\@psdoloop#1,#2,#3\@@#4#5{\def#4{#1}\ifx #4\@nnil \else
       #5\def#4{#2}\ifx #4\@nnil \else#5\@ipsdoloop #3\@@#4{#5}\fi\fi}
\def\@ipsdoloop#1,#2\@@#3#4{\def#3{#1}\ifx #3\@nnil
       \let\@nextwhile=\@psdonoop \else
      #4\relax\let\@nextwhile=\@ipsdoloop\fi\@nextwhile#2\@@#3{#4}}
\def\@tpsdo#1:=#2\do#3{\xdef\@psdotmp{#2}\ifx\@psdotmp\@empty \else
    \@tpsdoloop#2\@nil\@nil\@@#1{#3}\fi}
\def\@tpsdoloop#1#2\@@#3#4{\def#3{#1}\ifx #3\@nnil
       \let\@nextwhile=\@psdonoop \else
      #4\relax\let\@nextwhile=\@tpsdoloop\fi\@nextwhile#2\@@#3{#4}}
\ifx\undefined\fbox
\newdimen\fboxrule
\newdimen\fboxsep
\newdimen\ps@tempdima
\newbox\ps@tempboxa
\long\def\fbox#1{\leavevmode\setbox\ps@tempboxa\hbox{#1}\ps@tempdima\fboxrule
    \advance\ps@tempdima \fboxsep \advance\ps@tempdima \dp\ps@tempboxa
   \hbox{\lower \ps@tempdima\hbox
  {\vbox{\hrule height \fboxrule
          \hbox{\vrule width \fboxrule \hskip\fboxsep
          \vbox{\vskip\fboxsep \box\ps@tempboxa\vskip\fboxsep}\hskip
                 \fboxsep\vrule width \fboxrule}
                 \hrule height \fboxrule}}}}
\fi
\newread\ps@stream
\newif\ifnot@eof       
\newif\if@noisy        
\newif\if@atend        
\newif\if@psfile       
{\catcode`\%=12\global\gdef\epsf@start{
\def\epsf@PS{PS}
\def\epsf@getbb#1{%
\openin\ps@stream=#1
\ifeof\ps@stream\ps@typeout{Error, File #1 not found}\else
   {\not@eoftrue \chardef\other=12
    \def\do##1{\catcode`##1=\other}\dospecials \catcode`\ =10
    \catcode`\^^L=9 \catcode`\^^?=9
    \@psfiletrue
    \loop
       \if@psfile
	  \read\ps@stream to \epsf@fileline
       \else{
	  \obeyspaces
          \read\ps@stream to \epsf@tmp\global\let\epsf@fileline\epsf@tmp}
       \fi
       \ifeof\ps@stream\not@eoffalse\else
       \if@psfile\else
       \expandafter\epsf@test\epsf@fileline:. \\%
       \fi
          \expandafter\epsf@aux\epsf@fileline:. \\%
       \fi
   \ifnot@eof\repeat
   }\closein\ps@stream\fi}%
\long\def\epsf@test#1#2#3:#4\\{\def\epsf@testit{#1#2}
			\ifx\epsf@testit\epsf@start\else
\ps@typeout{Warning! File does not start with `\epsf@start'.  
It may not be a PostScript file.}
			\fi
			\@psfiletrue} 
{\catcode`\%=12\global\let\epsf@percent=
\long\def\epsf@aux#1#2:#3\\{\ifx#1\epsf@percent
   \def\epsf@testit{#2}\ifx\epsf@testit\epsf@bblit
	\@atendfalse
        \epsf@atend #3 . \\%
	\if@atend	
	   \if@verbose{
		\ps@typeout{psfig: found `(atend)'; continuing search}
	   }\fi
        \else
        \epsf@grab #3 . . . \\%
        \not@eoffalse
        \global\no@bbfalse
        \fi
   \fi\fi}%
\def\epsf@grab #1 #2 #3 #4 #5\\{%
   \global\def\epsf@llx{#1}\ifx\epsf@llx\empty
      \epsf@grab #2 #3 #4 #5 .\\\else
   \global\def\epsf@lly{#2}%
   \global\def\epsf@urx{#3}\global\def\epsf@ury{#4}\fi}%
\def\epsf@atendlit{(atend)}
\def\epsf@atend #1 #2 #3\\{%
   \def\epsf@tmp{#1}\ifx\epsf@tmp\empty
      \epsf@atend #2 #3 .\\\else
   \ifx\epsf@tmp\epsf@atendlit\@atendtrue\fi\fi}

\chardef\psletter = 11 
\chardef\other = 12

\newif\ifdebug 
\newif\ifc@mpute 
\c@mputetrue 

\let\then = \relax
\def\r@dian{pt }
\let\r@dians = \r@dian
\let\dimensionless@nit = \r@dian
\let\dimensionless@nits = \dimensionless@nit
\def\internal@nit{sp }
\let\internal@nits = \internal@nit
\newif\ifstillc@nverging
\def \Mess@ge #1{\ifdebug \then \message {#1} \fi}

{ 
	\catcode `\@ = \psletter
	\gdef \nodimen {\expandafter \n@dimen \the \dimen}
	\gdef \term #1 #2 #3%
	       {\edef \t@ {\the #1}
		\edef \t@@ {\expandafter \n@dimen \the #2\r@dian}%
		\t@rm {\t@} {\t@@} {#3}%
	       }
	\gdef \t@rm #1 #2 #3%
	       {{%
		\count 0 = 0
		\dimen 0 = 1 \dimensionless@nit
		\dimen 2 = #2\relax
		\Mess@ge {Calculating term #1 of \nodimen 2}%
		\loop
		\ifnum	\count 0 < #1
		\then	\advance \count 0 by 1
			\Mess@ge {Iteration \the \count 0 \space}%
			\Multiply \dimen 0 by {\dimen 2}%
			\Mess@ge {After multiplication, term = \nodimen 0}%
			\Divide \dimen 0 by {\count 0}%
			\Mess@ge {After division, term = \nodimen 0}%
		\repeat
		\Mess@ge {Final value for term #1 of
				\nodimen 2 \space is \nodimen 0}%
		\xdef \Term {#3 = \nodimen 0 \r@dians}%
		\aftergroup \Term
	       }}
	\catcode `\p = \other
	\catcode `\t = \other
	\gdef \n@dimen #1pt{#1} 
}

\def \Divide #1by #2{\divide #1 by #2} 

\def \Multiply #1by #2
       {{
	\count 0 = #1\relax
	\count 2 = #2\relax
	\count 4 = 65536
	\Mess@ge {Before scaling, count 0 = \the \count 0 \space and
			count 2 = \the \count 2}%
	\ifnum	\count 0 > 32767 
	\then	\divide \count 0 by 4
		\divide \count 4 by 4
	\else	\ifnum	\count 0 < -32767
		\then	\divide \count 0 by 4
			\divide \count 4 by 4
		\else
		\fi
	\fi
	\ifnum	\count 2 > 32767 
	\then	\divide \count 2 by 4
		\divide \count 4 by 4
	\else	\ifnum	\count 2 < -32767
		\then	\divide \count 2 by 4
			\divide \count 4 by 4
		\else
		\fi
	\fi
	\multiply \count 0 by \count 2
	\divide \count 0 by \count 4
	\xdef \product {#1 = \the \count 0 \internal@nits}%
	\aftergroup \product
       }}

\def\r@duce{\ifdim\dimen0 > 90\r@dian \then   
		\multiply\dimen0 by -1
		\advance\dimen0 by 180\r@dian
		\r@duce
	    \else \ifdim\dimen0 < -90\r@dian \then  
		\advance\dimen0 by 360\r@dian
		\r@duce
		\fi
	    \fi}

\def\Sine#1%
       {{%
	\dimen 0 = #1 \r@dian
	\r@duce
	\ifdim\dimen0 = -90\r@dian \then
	   \dimen4 = -1\r@dian
	   \c@mputefalse
	\fi
	\ifdim\dimen0 = 90\r@dian \then
	   \dimen4 = 1\r@dian
	   \c@mputefalse
	\fi
	\ifdim\dimen0 = 0\r@dian \then
	   \dimen4 = 0\r@dian
	   \c@mputefalse
	\fi
	\ifc@mpute \then
		\divide\dimen0 by 180
		\dimen0=3.141592654\dimen0
		\dimen 2 = 3.1415926535897963\r@dian 
		\divide\dimen 2 by 2 
		\Mess@ge {Sin: calculating Sin of \nodimen 0}%
		\count 0 = 1 
		\dimen 2 = 1 \r@dian 
		\dimen 4 = 0 \r@dian 
		\loop
			\ifnum	\dimen 2 = 0 
			\then	\stillc@nvergingfalse
			\else	\stillc@nvergingtrue
			\fi
			\ifstillc@nverging 
			\then	\term {\count 0} {\dimen 0} {\dimen 2}%
				\advance \count 0 by 2
				\count 2 = \count 0
				\divide \count 2 by 2
				\ifodd	\count 2 
				\then	\advance \dimen 4 by \dimen 2
				\else	\advance \dimen 4 by -\dimen 2
				\fi
		\repeat
	\fi		
			\xdef \sine {\nodimen 4}%
       }}

\def\Cosine#1{\ifx\sine\UnDefined\edef\Savesine{\relax}\else
		             \edef\Savesine{\sine}\fi
	{\dimen0=#1\r@dian\advance\dimen0 by 90\r@dian
	 \Sine{\nodimen 0}
	 \xdef\cosine{\sine}
	 \xdef\sine{\Savesine}}}

\def\psdraft{
	\def\@psdraft{0}
}
\def\psfull{
	\def\@psdraft{100}
}

\psfull

\newif\if@scalefirst
\def\psscalefirst{\@scalefirsttrue}
\def\psrotatefirst{\@scalefirstfalse}
\psrotatefirst

\newif\if@draftbox
\def\psnodraftbox{
	\@draftboxfalse
}
\def\psdraftbox{
	\@draftboxtrue
}
\@draftboxtrue

\newif\if@prologfile
\newif\if@postlogfile
\def\pssilent{
	\@noisyfalse
}
\def\psnoisy{
	\@noisytrue
}
\psnoisy
\newif\if@bbllx
\newif\if@bblly
\newif\if@bburx
\newif\if@bbury
\newif\if@height
\newif\if@width
\newif\if@rheight
\newif\if@rwidth
\newif\if@angle
\newif\if@clip
\newif\if@verbose
\def\@p@@sclip#1{\@cliptrue}

\newif\if@decmpr

\def\@p@@sfigure#1{\def\@p@sfile{null}\def\@p@sbbfile{null}
	        \openin1=#1.bb
		\ifeof1\closein1
	        	\openin1=\figurepath#1.bb
			\ifeof1\closein1
			        \openin1=#1
				\ifeof1\closein1%
				       \openin1=\figurepath#1
					\ifeof1
					   \ps@typeout{Error, File #1 not found}
						\if@bbllx\if@bblly
				   		\if@bburx\if@bbury
			      				\def\@p@sfile{#1}%
			      				\def\@p@sbbfile{#1}%
							\@decmprfalse
				  	   	\fi\fi\fi\fi
					\else\closein1
				    		\def\@p@sfile{\figurepath#1}%
				    		\def\@p@sbbfile{\figurepath#1}%
						\@decmprfalse
	                       		\fi%
			 	\else\closein1%
					\def\@p@sfile{#1}
					\def\@p@sbbfile{#1}
					\@decmprfalse
			 	\fi
			\else
				\def\@p@sfile{\figurepath#1}
				\def\@p@sbbfile{\figurepath#1.bb}
				\@decmprtrue
			\fi
		\else
			\def\@p@sfile{#1}
			\def\@p@sbbfile{#1.bb}
			\@decmprtrue
		\fi}

\def\@p@@sfile#1{\@p@@sfigure{#1}}

\def\@p@@sbbllx#1{
		\ps@typeout{bbllx is #1} 
		\@bbllxtrue
		\dimen100=#1
		\edef\@p@sbbllx{\number\dimen100}
}
\def\@p@@sbblly#1{
		\ps@typeout{bblly is #1} 
		\@bbllytrue
		\dimen100=#1
		\edef\@p@sbblly{\number\dimen100}
}
\def\@p@@sbburx#1{
		\ps@typeout{bburx is #1} 
		\@bburxtrue
		\dimen100=#1
		\edef\@p@sbburx{\number\dimen100}
}
\def\@p@@sbbury#1{
		\ps@typeout{bbury is #1} 
		\@bburytrue
		\dimen100=#1
		\edef\@p@sbbury{\number\dimen100}
}
\def\@p@@sheight#1{
		\@heighttrue
		\dimen100=#1
   		\edef\@p@sheight{\number\dimen100}
		\ps@typeout{Height is \@p@sheight} 
}
\def\@p@@swidth#1{
		\ps@typeout{Width is #1} 
		\@widthtrue
		\dimen100=#1
		\edef\@p@swidth{\number\dimen100}
}
\def\@p@@srheight#1{
		\ps@typeout{Reserved height is #1} 
		\@rheighttrue
		\dimen100=#1
		\edef\@p@srheight{\number\dimen100}
}
\def\@p@@srwidth#1{
		\ps@typeout{Reserved width is #1} 
		\@rwidthtrue
		\dimen100=#1
		\edef\@p@srwidth{\number\dimen100}
}
\def\@p@@sangle#1{
		\ps@typeout{Rotation is #1} 
		\@angletrue
		\edef\@p@sangle{#1} 
}
\def\@p@@ssilent#1{
		\@verbosefalse
}
\def\@p@@sprolog#1{\@prologfiletrue\def\@prologfileval{#1}}
\def\@p@@spostlog#1{\@postlogfiletrue\def\@postlogfileval{#1}}
\def\@cs@name#1{\csname #1\endcsname}
\def\@setparms#1=#2,{\@cs@name{@p@@s#1}{#2}}
\def\ps@init@parms{
		\@bbllxfalse \@bbllyfalse
		\@bburxfalse \@bburyfalse
		\@heightfalse \@widthfalse
		\@rheightfalse \@rwidthfalse
		\def\@p@sbbllx{}\def\@p@sbblly{}
		\def\@p@sbburx{}\def\@p@sbbury{}
		\def\@p@sheight{}\def\@p@swidth{}
		\def\@p@srheight{}\def\@p@srwidth{}
		\def\@p@sangle{0}
		\def\@p@sfile{} \def\@p@sbbfile{}
		\def\@p@scost{10}
		\def\@sc{}
		\@prologfilefalse
		\@postlogfilefalse
		\@clipfalse
		\if@noisy
			\@verbosetrue
		\else
			\@verbosefalse
		\fi
}
\def\parse@ps@parms#1{
	 	\@psdo\@psfiga:=#1\do
		   {\expandafter\@setparms\@psfiga,}}
\newif\ifno@bb
\def\bb@missing{
	\if@verbose{
		\ps@typeout{psfig: searching \@p@sbbfile \space  for bounding box}
	}\fi
	\no@bbtrue
	\epsf@getbb{\@p@sbbfile}
        \ifno@bb \else \bb@cull\epsf@llx\epsf@lly\epsf@urx\epsf@ury\fi
}	
\def\bb@cull#1#2#3#4{
	\dimen100=#1 bp\edef\@p@sbbllx{\number\dimen100}
	\dimen100=#2 bp\edef\@p@sbblly{\number\dimen100}
	\dimen100=#3 bp\edef\@p@sbburx{\number\dimen100}
	\dimen100=#4 bp\edef\@p@sbbury{\number\dimen100}
	\no@bbfalse
}
\newdimen\p@intvaluex
\newdimen\p@intvaluey
\def\rotate@#1#2{{\dimen0=#1 sp\dimen1=#2 sp
		  \global\p@intvaluex=\cosine\dimen0
		  \dimen3=\sine\dimen1
		  \global\advance\p@intvaluex by -\dimen3
		  \global\p@intvaluey=\sine\dimen0
		  \dimen3=\cosine\dimen1
		  \global\advance\p@intvaluey by \dimen3
		  }}
\def\compute@bb{
		\no@bbfalse
		\if@bbllx \else \no@bbtrue \fi
		\if@bblly \else \no@bbtrue \fi
		\if@bburx \else \no@bbtrue \fi
		\if@bbury \else \no@bbtrue \fi
		\ifno@bb \bb@missing \fi
		\ifno@bb \ps@typeout{FATAL ERROR: no bb supplied or found}
			\no-bb-error
		\fi
\ps@typeout{BB: \@p@sbbllx, \@p@sbblly, \@p@sbburx, \@p@sbbury} 
		\count203=\@p@sbburx
		\count204=\@p@sbbury
		\advance\count203 by -\@p@sbbllx
		\advance\count204 by -\@p@sbblly
		\edef\ps@bbw{\number\count203}
		\edef\ps@bbh{\number\count204}
		\if@angle
			\Sine{\@p@sangle}\Cosine{\@p@sangle}
	        	{\dimen100=\maxdimen\xdef\r@p@sbbllx{\number\dimen100}
					    \xdef\r@p@sbblly{\number\dimen100}
			                    \xdef\r@p@sbburx{-\number\dimen100}
					    \xdef\r@p@sbbury{-\number\dimen100}}
                        \def\minmaxtest{
			   \ifnum\number\p@intvaluex<\r@p@sbbllx
			      \xdef\r@p@sbbllx{\number\p@intvaluex}\fi
			   \ifnum\number\p@intvaluex>\r@p@sbburx
			      \xdef\r@p@sbburx{\number\p@intvaluex}\fi
			   \ifnum\number\p@intvaluey<\r@p@sbblly
			      \xdef\r@p@sbblly{\number\p@intvaluey}\fi
			   \ifnum\number\p@intvaluey>\r@p@sbbury
			      \xdef\r@p@sbbury{\number\p@intvaluey}\fi
			   }
			\rotate@{\@p@sbbllx}{\@p@sbblly}
			\minmaxtest
			\rotate@{\@p@sbbllx}{\@p@sbbury}
			\minmaxtest
			\rotate@{\@p@sbburx}{\@p@sbblly}
			\minmaxtest
			\rotate@{\@p@sbburx}{\@p@sbbury}
			\minmaxtest
			\edef\@p@sbbllx{\r@p@sbbllx}\edef\@p@sbblly{\r@p@sbblly}
			\edef\@p@sbburx{\r@p@sbburx}\edef\@p@sbbury{\r@p@sbbury}
		\fi
		\count203=\@p@sbburx
		\count204=\@p@sbbury
		\advance\count203 by -\@p@sbbllx
		\advance\count204 by -\@p@sbblly
		\edef\@bbw{\number\count203}
		\edef\@bbh{\number\count204}
}
\def\in@hundreds#1#2#3{\count240=#2 \count241=#3
		     \count100=\count240	
		     \divide\count100 by \count241
		     \count101=\count100
		     \multiply\count101 by \count241
		     \advance\count240 by -\count101
		     \multiply\count240 by 10
		     \count101=\count240	
		     \divide\count101 by \count241
		     \count102=\count101
		     \multiply\count102 by \count241
		     \advance\count240 by -\count102
		     \multiply\count240 by 10
		     \count102=\count240	
		     \divide\count102 by \count241
		     \count200=#1\count205=0
		     \count201=\count200
			\multiply\count201 by \count100
		 	\advance\count205 by \count201
		     \count201=\count200
			\divide\count201 by 10
			\multiply\count201 by \count101
			\advance\count205 by \count201
		     \count201=\count200
			\divide\count201 by 100
			\multiply\count201 by \count102
			\advance\count205 by \count201
		     \edef\@result{\number\count205}
}
\def\compute@wfromh{
		\in@hundreds{\@p@sheight}{\@bbw}{\@bbh}
		\edef\@p@swidth{\@result}
}
\def\compute@hfromw{
	        \in@hundreds{\@p@swidth}{\@bbh}{\@bbw}
		\edef\@p@sheight{\@result}
}
\def\compute@handw{
		\if@height
			\if@width
			\else
				\compute@wfromh
			\fi
		\else
			\if@width
				\compute@hfromw
			\else
				\edef\@p@sheight{\@bbh}
				\edef\@p@swidth{\@bbw}
			\fi
		\fi
}
\def\compute@resv{
		\if@rheight \else \edef\@p@srheight{\@p@sheight} \fi
		\if@rwidth \else \edef\@p@srwidth{\@p@swidth} \fi
}
\def\compute@sizes{
	\compute@bb
	\if@scalefirst\if@angle
	\if@width
	   \in@hundreds{\@p@swidth}{\@bbw}{\ps@bbw}
	   \edef\@p@swidth{\@result}
	\fi
	\if@height
	   \in@hundreds{\@p@sheight}{\@bbh}{\ps@bbh}
	   \edef\@p@sheight{\@result}
	\fi
	\fi\fi
	\compute@handw
	\compute@resv}

\def\psfig#1{\vbox {
	\ps@init@parms
	\parse@ps@parms{#1}
	\compute@sizes
	\ifnum\@p@scost<\@psdraft{
		\special{ps::[begin] 	\@p@swidth \space \@p@sheight \space
				\@p@sbbllx \space \@p@sbblly \space
				\@p@sbburx \space \@p@sbbury \space
				startTexFig \space }
		\if@angle
			\special {ps:: \@p@sangle \space rotate \space}
		\fi
		\if@clip{
			\if@verbose{
				\ps@typeout{(clip)}
			}\fi
			\special{ps:: doclip \space }
		}\fi
		\if@prologfile
		    \special{ps: plotfile \@prologfileval \space } \fi
		\if@decmpr{
			\if@verbose{
				\ps@typeout{psfig: including \@p@sfile.Z \space }
			}\fi
			\special{ps: plotfile "`zcat \@p@sfile.Z" \space }
		}\else{
			\if@verbose{
				\ps@typeout{psfig: including \@p@sfile \space }
			}\fi
			\special{ps: plotfile \@p@sfile \space }
		}\fi
		\if@postlogfile
		    \special{ps: plotfile \@postlogfileval \space } \fi
		\special{ps::[end] endTexFig \space }
		\vbox to \@p@srheight sp{
			\hbox to \@p@srwidth sp{
				\hss
			}
		\vss
		}
	}\else{
		\if@draftbox{		
			\hbox{\frame{\vbox to \@p@srheight sp{
			\vss
			\hbox to \@p@srwidth sp{ \hss \@p@sfile \hss }
			\vss
			}}}
		}\else{
			\vbox to \@p@srheight sp{
			\vss
			\hbox to \@p@srwidth sp{\hss}
			\vss
			}
		}\fi

	}\fi
}}
\psfigRestoreAt
\let\@=\LaTeXAtSign